\documentclass[12pt]{article} 
\usepackage{amsmath,amssymb,eucal} 
\usepackage{bbm}  
\font\tenrm=cmr10

\font\cmssl=cmss10 at 12 pt  
   
\font\bigss=cmssdc10 scaled 2300

\font\cmsslll=cmss10 at 14 pt

\renewcommand{\a}{\alpha}

\renewcommand{\d}{\delta}

\renewcommand{\o}{\omega}

\newcommand{\s}{\sigma}

\newcommand{\G}{\Gamma}

\newcommand{\bR}{\mathbb{R}}

\newcommand{\id}   {{\mathbbm{1}}}

\renewcommand{\square}{\kern1pt\vbox  
               {\hrule height 0.6pt\hbox{\vrule width 0.6pt\hskip 3pt  
    \vbox{\vskip 6pt}\hskip 3pt\vrule width 0.6pt}\hrule height0.6pt}  
    \kern1pt}

\newcommand{\ra}{\rightarrow}

\newtheorem{Th}{Theorem}  
\newtheorem{Prop}{Proposition}  
\newtheorem{Cor}{Corollary}  
\newtheorem{Lem}{Lemma}  
\newtheorem{Def}{Definition}  
\newtheorem{Ex}{Example}  
\newcommand{\bt}{\begin{Th}\ \ }  
\newcommand{\et}{\end{Th}}  
\newcommand{\bp}{\begin{Prop}\ \ }  
\newcommand{\ep}{\end{Prop}}  
\newcommand{\bc}{\begin{Cor}\ \ }  
\newcommand{\ec}{\end{Cor}}  
\newcommand{\bl}{\begin{Lem}\ \ }  
\newcommand{\el}{\end{Lem}}  
\newcommand{\bd}{\begin{Def}\ \ }  
\newcommand{\ed}{\end{Def}}  
  
\newcommand{\pf}{\noindent{\it Proof:\ \ }}  
\newcommand{\qed}{\hfill\square}  
\newcommand{\n}{\nabla}

\newcommand{\be}{\begin{equation}}  
\newcommand{\ee}{\end{equation}}  
  
\newcommand\re[1]{(\ref{#1})}  
\newcommand{\arr}{\begin{array}{rlll}}  
\newcommand{\ea}{\end{array}}  
\newcommand{\bea}{\begin{eqnarray}}  
\newcommand{\eea}{\end{eqnarray}}  
\newcommand{\bean}{\begin{eqnarray*}}  
\newcommand{\eean}{\end{eqnarray*}}  
  
\catcode`@=11  
\@addtoreset{equation}{section}  
\catcode`@=12

\begin{document}  
\begin{titlepage}
 \rightline{} 
\vskip 1.5 true cm  
\begin{center}  
{\bigss  Geometric construction of the r-map:  
from affine\\[.5em] special real to special K\"ahler manifolds}  
\vskip 1.0 true cm   
{\cmsslll  D.V.\ Alekseevsky and  V.\ Cort\'es} \\[3pt] 
{\tenrm   The University of Edinburgh and Maxwell Institute for
Mathematical Sciences\\ 
JCMB, The King's buildings,  
Edinburgh, EH9 3JZ, UK \\
D.Aleksee@ed.ac.uk}\\[1em]  
{\tenrm   Department Mathematik 
und Zentrum f\"ur Mathematische Physik\\ 
Universit\"at Hamburg, 
Bundesstra{\ss}e 55, 
D-20146 Hamburg, Germany\\  
cortes@math.uni-hamburg.de}\\[1em]   
October 28, 2008 
\end{center}  
\vskip 1.0 true cm  
\baselineskip=18pt  
\begin{abstract}  
\noindent  
We give an intrinsic definition of (affine very) special real manifolds   
and realise any such manifold $M$ as a domain in affine 
space equipped with a metric which is the Hessian of a cubic polynomial. 
We prove that the tangent bundle $N=TM$ carries a canonical structure  
of (affine) special K\"ahler manifold. This gives an intrinsic 
description of the $r$-map as the map 
$M\mapsto N=TM$. On the physics side, 
this map corresponds to the dimensional reduction of rigid vector 
multiplets from 5 to 4 space-time dimensions. 
We generalise this construction to the case when $M$ is any Hessian
manifold.  
\end{abstract}

\end{titlepage}  
\tableofcontents
\section*{Introduction} 
{\it Projective (or local) very special real geometry} is the  
scalar geometry of five-dimensional supergravity coupled to vector multiplets 
\cite{GST,dWvP,ACDV}.  We will usually omit the word ``very''. 
It is locally the geometry of a nondegenerate 
hypersurface ${\cal H}\subset \bR^{n+1}$ defined by a homogeneous cubic
polynomial $h$. In relation with string compactifications 
the polynomial $h$ could be, for instance, the cubic form 
$$ h([\a ]) = \int_X\a \wedge \a \wedge \a$$
on the second cohomology of a Calabi-Yau 3-fold $X$.  

In this paper we are concerned with 
{\it affine (or rigid) very special real geometry}, which
is the scalar geometry of five-dimensional rigid 
vector multiplets \cite{CMMS1}. The Lagrangian of rigid 
vector multiplets is encoded in a (not necessarily homogeneous) 
cubic polynomial and the metric of the scalar manifold is the Hessian of this 
polynomial.  

In the first part of the paper we shall 
provide an intrinsic definition of the notion of an 
affine special real manifold and study, in particular, 
the geometry of such manifolds:
\bd 
{\it An affine special real manifold} $(M,g,\n )$ is  
a pseudo-Riemannian manifold $(M,g)$ endowed with a flat torsion-free
connection $\n$ such that the tensor
field $S=\n g$ is  totally symmetric and $\n$-parallel. 
\ed
We relate the intrinsic definition to the description in the physical 
literature in terms of a 
cubic prepotential. In fact, we show  that any 
simply connected affine special real manifold $(M,g,\n )$ of dimension 
$n$ admits an affine immersion 
$\psi$ onto a domain $V\subset \bR^n$, such that $g$ is the pull back of 
the Hessian
of a cubic polynomial $h$, see Corollary \ref{immCor}. The pair $(V,h)$ is 
unique up to affine tranformations of $\bR^n$.  We calculate the curvature
tensor of a special real manifold (and, more generally, of a Hessian manifold, 
see Corollary \ref{flatCor}) and find a simple expression in terms of  
the tensor $S$. As an application, we obtain 
that a special real manifold with a 
positive definite metric has 
nonnegative Ricci curvature, see Corollary \ref{RicCor}. 

Dimensional reduction of (local/rigid) supersymmetric theories from 5 to 4 
space-time dimensions induces a correspondence between the 
respective scalar geometries, which is know as the (local/rigid) 
{\it r-map} \cite{dWvP,CMMS1}. The relevant scalar geometry
in 4 space-time dimensions is (projective/affine) 
{\it special K\"ahler geometry}, see \cite{C} for a survey.
The (local/rigid) r-map associates a (projective/affine) special K\"ahler
manifold to any (projective/affine) special real  manifold. 
The r-map is explicitely known in terms of the 
prepotentials, which locally define special real and special K\"ahler geometry.
In the affine case, for instance, the r-map associates to the real 
cubic polynomial $h(x^1,\ldots ,x^n)$ defining the special 
real manifold the holomorphic prepotential  
$F(z^1,\ldots ,z^n)= \frac{1}{2i}h(z^1,\ldots ,z^n)$ 
defining the corresponding special K\"ahler manifold in terms of special
holomorphic coordinates $z^1,\ldots, z^n$ 
\cite{CMMS1}. 
However, an intrinsic geometric construction of the affine and projective 
r-maps is missing. In the second part of the paper we shall give such 
a construction in the affine 
case\footnote{The projective case is the subject of work in progress.}.

We show that the tangent bundle $N=TM$ of any (affine) 
special real manifold  $(M,g,\n )$ 
carries the structure of an (affine) 
special K\"ahler manifold $(N,J,g^N,\n^N)$.  

More precisely,  the 
special K\"ahler structure $(J,g^N,\n^N)$ on $N$ is canonically
associated to the geometric data $(g,\n )$ on $M$, see Theorem \ref{rmapThm}. 
Recall that a special K\"ahler manifold  $(N,J,g^N,\n^N)$ is 
(pseudo-)K\"ahler manifold   $(N,J,g^N)$ endowed with a connection
$\n^N$ which is {\it special} (i.e.\ $\n^NJ$ is symmetric), 
torsion-free,  symplectic 
(with respect to the K\"ahler form) and flat. 
The map \begin{eqnarray}\label{rmapEqu} 
 \mathbf{r} : \{ \mbox{special real manifolds}\} &\ra& 
\{ \mbox{special K\"ahler manifolds}\}\\
 (M,g,\n )&\mapsto& (N,J,g^N,\n^N),\nonumber
\end{eqnarray}
which associates to the special real manifold $(M,g,\n )$
the special K\"ahler manifold\linebreak[4] 
$(N=TM,J,g^N,\n^N)$ is our geometric
construction of the r-map. As an application, we prove that there is no compact
simply connected special real manifold with a positive definite metric, see
Theorem \ref{lastThm}.

We show that our r-map extends naturally to a map 
\begin{eqnarray}\label{extrmapEqu}
 \mathbf{r} : \{ \mbox{Hessian manifolds}\} &\ra&
\left\{ \begin{array}{c}\mbox{K\"ahler manifolds  with a torsion-free,}\\
\mbox{symplectic and special connection}\end{array}\right\}\\
(M,g,\n )&\mapsto& (N,J,g^N,\n^N),\nonumber
\end{eqnarray}
A {\it Hessian manifold} $(M,g,\n)$ is a pseudo-Riemannian manifold such that
$S=\n g$ is totally symmetric (but not necessarily $\n$-parallel). 
The flatness of the connection $\n^N$ is lost
when the r-map is applied to Hessian manifolds which are not
special real. In fact, $\n^N$ is flat if and only if $(M,g,\n)$ is 
special real. Moreover, the manifold $(N,J,g^N,\n^N)$ associated
to a Hessian manifold $(M,g,\n)$ by the r-map is again Hessian 
if and only if  $(M,g,\n)$ is special real, see Corollary \ref{8.5Cor}. 

Finally, we characterise the image of the
maps \re{rmapEqu}  and \re{extrmapEqu} in Theorems \ref{characThm} and 
\ref{characconvThm}. We calculate the curvature tensors of the
Levi-Civita connection of the metric 
$g_N$ and of $\n^N$, which have a simple expression in terms of the symmetric 
tensor $S=\n g$ and $\n S$, respectively, see Corollaries \ref{curvCor},
\ref{Cor9} and \ref{Cor10}.   

In particular, it follows from those theorems that 
a special K\"ahler manifold\linebreak[4] 
$(N,J,g^N,\n^N)$ of real dimension $2n$ can be locally obtained from the 
r-map if and only if it admits locally $n$  
holomorphic\footnote{Recall that a real 
vector field $X$ on a complex manifold $(N,J)$ is called {\it holomorphic} if 
${\cal L}_XJ=0$.} Killing vector fields which span a Lagrangian
distribution and which are $\n^N$-parallel along this distribution. 

We prove that under some assumptions a simply connected
$n$-dimensional Hessian manifold admits a canonical realisation as an 
improper affine hypersphere in $\bR^{n+1}$ equipped with the Blaschke
metric and the induced connection.

\section{Hessian geometry and affine special real geometry}
\bd \label{mainDef} A {\cmssl Hessian manifold} (cf.\ \cite{S}) 
$(M,g,\n )$ is a pseudo-Riemannian manifold 
$(M,g)$ with a flat torsion-free connection $\n$ such that $S=\n g$ is 
a symmetric three-form (cubic form). An {\cmssl affine (very) 
special real manifold}  $ (M,g,\n )$ 
is a Hessian manifold $(M,g,\n )$ with $\n$-parallel cubic form $S$. 
\ed 

\begin{Ex} \label{1stEx}
Let $h$ be a smooth function in affine space $V\cong \bR^n$. We will
say that $h$ is {\cmssl nondegenerate} at a point $x_0\in V$ if the
Hessian $\partial^2h(x_0)$ is nondegenerate, where $\partial$ 
is the flat connection
in $V$. We denote by $V(x_0) \subset V$ the connected component of
$x_0$ in $\{ \det \partial^2h\neq 0\}\subset V$. The domain
$V(x_0)$ is equipped with the pseudo-Riemannian metric 
$g=\partial^2 h$. Then $(V(x_0),g,\partial )$ is 
a Hessian manifold. Indeed $S=\partial g=\partial^3h$ is completely
symmetric. It is an 
affine special real manifold if and only if the cubic form $S$ is constant.
\end{Ex} 
\bp \label{PropAffSpecMan} Any Hessian manifold (respectively, 
affine special real manifold) $(M,g,\n )$ is locally
isomorphic to a domain $(V(x_0),g,\partial )$ associated with a 
smooth function $h$ (respectively, 
cubic polynomial $h$), as in Example \ref{1stEx}.  The polynomial is 
given by 
$$h= \frac{1}{6}\sum S_{ijk}x^ix^jx^k + \frac{1}{2}b_{ij}x^ix^j.$$  
Here $x^i$ are $\n$-affine coordinates $g=\sum g_{ij}dx^idx^j$, 
$g_{ij}= \sum S_{ijk}x^k + b_{ij}$ and\linebreak  
$S= \sum S_{ijk}dx^idx^jdx^k$. For any Hessian manifold the 
n linearly independent gradient vector fields ${\rm grad}(x^i)$ commute
and the coordinate vector fields $\frac{\partial}{\partial x^i}$ are 
also commuting gradient vector fields. 
\ep 

\pf Since $S=\nabla g$ is totally symmetric, there exists locally 
a smooth function $h$ such that $g=\n^2h$. Moreover, if $\n S=0$, 
then the function $h$ is 
a cubic  polynomial in affine local coordinates $x^i$. 
Then $\partial^2h=\n^2 h=g=\sum g_{ij}dx^idx^j$ and $\partial^3h=\nabla^3h=S=
\sum S_{ijk}dx^idx^jdx^k$, where $g_{ij}=\sum a_{ijk}x^k+b_{ij}$ and 
$S_{ijk}=a_{ijk}$. This shows
that $h$ coincides with the above expression up to terms of degree less or
equal to 1 in the coordinates $x^i$, which do not contribute to
$g$ and $S$. Now we check that the vector fields ${\rm grad}(x^i) = 
\sum g^{ij}\partial_j$ commute: 
\begin{eqnarray}\label{symmEqu} [g^{ij}\partial_j,g^{kl}\partial_l] &=& 
g^{ij}\partial_jg^{kl}\partial_l
- g^{kl}\partial_lg^{ij}\partial_j = -g^{ij}g^{kk'}S_{k'l'j}g^{l'l}\partial_l
+g^{kl}g^{ii'}S_{i'j'l}g^{j'j}\partial_j\nonumber \\ 
&=& -S^{kli}\partial_l+S^{ilk}\partial_l
=0,
\end{eqnarray}
by the complete symmetry of $S$. (Here and in later 
calculations we use the Einstein summation 
convention.) The coordinate vector field $\frac{\partial}{\partial x^i}$ 
is the gradient of the function $\frac{\partial h}{\partial x^i}$. 
\qed

\bc \label{immCor} Let $(M,g,\n )$ be a simply connected 
Hessian manifold of dimension n. 
Then there exists an affine immersion 
$\psi :(M,\n ) \ra (\bR^n,\partial )$ onto
some domain $\psi (M) \subset \bR^n$,  
unique up to affine transformations of $\bR^n$. 
The gradients  ${\rm grad}(x^i)$ of the coordinate functions 
$x^i=\psi^i$ span a canonical n-dimensional 
commutative Lie algebra of vector 
fields. Conversely, a pseudo-Riemannian manifold
$(M,g)$ with n pointwise linearly independent gradient 
vector fields ${\rm grad}(x^i)$ is canonically extended
to a Hessian manifold $(M,g,\n )$.  

If $(M,g,\n )$ is an affine special real manifold, then, moreover, 
there exists a unique cubic polynomial $h$ on $\bR^n$ 
without linear and constant terms such that 
$g=\psi^*\partial^2h$.  
\ec 

\pf Since $M$ is simply connected, there exists a $\n$-parallel
coframe $(\xi^1, \cdots, \xi^n)$. The one-forms $\xi^i$ are
closed and, hence, exact, i.e.\ $\xi^i=dx^i$ for globally
defined functions $x^i$. These functions define an affine immersion
$\psi =(x^1,\cdots, x^n)$ since $\n dx^i=0$. 
Given a pseudo-Riemannian manifold
$(M,g)$ with n pointwise linearly independent gradient 
vector fields ${\rm grad}(x^i)$, there is a unique flat
torsion-free connection $\n$ such that $\n dx^i=0$. Then
$\n g$ is completely symmetric by \re{symmEqu}.

On any 
domain $U\subset M$ such that $\psi|_U : U \ra \psi(U)$ is a diffeomorphism 
there exists a smooth function $h_U$ such that 
$g_{ij}|_U= \partial_i\partial_jh_U$. The function $h_U$ is unique
up to the addition of an affine function in the local coordinates $x^i$. 
In the special real case 
$h_U$ is a cubic polynomial, which can be canonically chosen by 
the requirement that the linear and constant terms vanish: 
$h_U=\frac{1}{6}a_{ijk}x^ix^jx^k + \frac{1}{2}b_{ij}x^ix^j.$
The coefficients $a_{ijk}$, $b_{ij}$ are independent of $U$, since $M$ 
is connected and $h_U=h_V$ on overlaps $U\cap V\neq\emptyset$.    
Therefore $h=\frac{1}{6}a_{ijk}x^ix^jx^k + \frac{1}{2}b_{ij}x^ix^j$
is canonically associated to the immersion $\psi$ and 
satisfies $g=\psi^*\partial^2h$.  
\qed 

\noindent 
{\bf Remark:} Shima and Yagi proved that the domain $\psi (M)$ is
convex if $g$ is positive definite, see \cite{SY}. 

\bc A pseudo-Riemannian manifold $(M,g)$ admits the structure
of an affine special real manifold if it admits an atlas with 
affine transition functions such that the metric 
coefficients $g_{ij}=a_{ijk}x^k+ b_{ij}$ are affine functions and
the coefficients $a_{ijk}$ are symmetric. Then $g_{ij}= \partial_i
 \partial_j h $ where $h = \frac{1}{6}a_{ijk}x^i x^j x^k + \frac{1}{2}b_{ij}x^i x^j$
\ec 

\bc Let $(M,g,\n )$ be a simply connected Hessian manifold with transitive
action of $G={\rm Aut} (M,g,\n )$.  Then the affine immersion
$\psi : M \ra \psi (M) \subset \bR^n$ of Corollary \ref{immCor}
is a $G$-equivariant covering $M=G/H \ra \psi (M)=Gx_0=G/G_{x_0}$
over the (open) orbit $Gx_0$ of
a point $x_0\in \psi (M)$ with respect to an affine action $\a$ of
$G$ on $\bR^{n}$. In particular, there is no nonconstant $G$-invariant
function on the domain $\psi (M)$ and  at most one  (up to
scaling) relative invariant 
with character $\chi (a) = (\det A)^{-2}$, where $\a (a)x = Ax +b$, $a\in G$.   
If such a relative invariant $\d : \psi (M) \ra \bR$ exists, it is given 
(up 
to a constant factor)  by 
the formula $\delta \circ \psi 
= \det (g_{ij})$ (which in general defines only a local 
relative invariant). For an affine special real manifold, 
the (globally defined) relative invariant $\delta$ is a polynomial 
of degree $n$. 
\ec 

Recall that given a pseudo-Riemannian metric $g$ and a connection
$\n$ on manifold $M$, the {\cmssl conjugate connection}
$\bar{\n}$ is defined by
\[ Xg(Y,Z) =g(\n_XY,Z) + g(Y,\bar{\n}_XZ),\]
where $X, Y, Z$ are vector fields on $M$. 
Notice that $\bar{\n}_X=D_X+\hat{S}_X^*$, where $D$ is the Levi-Civita
connection, $\hat{S}=D-\n$ and
$\hat{S}_X^*$ is the metric adjoint of the endomorphism $\hat{S}_X$. 
\bp \label{conjProp} Let $(M,g,\n )$ be a Hessian manifold with cubic form 
$S=\n g$. Then $\hat{S}_X=\frac{1}{2}g^{-1}\circ S_X=\hat{S}_X^*$. 
The conjugate connection  is flat
and torsion-free and we have the following formulas: 
\begin{eqnarray*} \bar{\n}_X&=&D_X+\hat{S}_X\\
\n_X&=&D_X-\hat{S}_X. 
\end{eqnarray*}
\ep 

\pf 
The connection $\n + \frac{1}{2}g^{-1}\circ S$ is torsion-free, by the
symmetry of $S$. We check that it preserves the metric $g$:
$$ \nabla_X g +\frac{1}{2}(g^{-1}\circ S_X)\cdot g = S_X-
 \frac{1}{2} g(g^{-1}\circ S_X \cdot, \cdot) -\frac{1}{2} g(\cdot,g^{-1}\circ 
S_X \cdot)=
S_X-S_X=0.
$$
This shows that $\n + \frac{1}{2}g^{-1}\circ S$ is the Levi-Civita 
connection $D=\n +\hat{S}$. Hence 
$\frac{1}{2}g^{-1}\circ S_X=\hat{S}_X=\hat{S}_X^*$. 
It is clear that the conjugate connection $\bar{\n}_X=D_X+\hat{S}_X$
has zero torsion. It remains to calculate its curvature $\bar{R}$. 
For this we write $\bar{\n}=\n +2\hat{S}$ and compute:
\begin{eqnarray} \label{nablaSEqu} \bar{R}(X,Y)&=& 
R^{\nabla}(X,Y)+2d^{\nabla}\hat{S}(X,Y)+4[\hat{S}_X,
\hat{S}_Y]\nonumber \\
&=&2(d^{\nabla}\hat{S}(X,Y)+2[\hat{S}_X,
\hat{S}_Y])\nonumber\\
\nabla_X\hat{S} &=&\frac{1}{2}\nabla_X (g^{-1}\circ S)=-\frac{1}{2}g^{-1}\circ S_Xg^{-1}\circ S + \frac{1}{2}g^{-1}\circ \n_XS\nonumber\\
&=&
-2\hat{S}_X\hat{S}+ \frac{1}{2}g^{-1}\circ \n_XS. 
\end{eqnarray} 
Therefore 
\begin{eqnarray*}d^{\nabla}\hat{S}(X,Y)&=&(\nabla_X\hat{S})Y-(\nabla_Y\hat{S})X=-2
[\hat{S}_X\hat{S}_Y],
\end{eqnarray*} 
since $\n S$ is symmetric. Thus $\bar{R}=0$. \qed  

\bc \label{flatCor} 
Under the assumptions of the proposition, the following formulas
are satisfied:
\begin{eqnarray*}
R^D(X,Y)&=&-[\hat{S}_X,\hat{S}_Y],\quad 
d^D\hat{S}=0,\quad 
d^\n \hat{S}+2[\hat{S},\hat{S}]=0.
\end{eqnarray*}  
\ec 

\pf 
The first two equations are obtained by taking the sum and difference of
the equations
\begin{eqnarray*}
0&=&\bar{R}=R^D+d^D\hat{S}+ [\hat{S},\hat{S}]\\
0&=&R^\n =R^D-d^D\hat{S} + [\hat{S},\hat{S}]. 
\end{eqnarray*} 
The third equation follows from $\bar{\n}=\n + 2\hat{S}$:
$$0=\bar{R}=R^\n+2d^\n \hat{S}+4[\hat{S},\hat{S}]=2(d^\n \hat{S}+2[\hat{S},\hat{S}]).$$
\qed

The following results are analogues of the corresponding
results in special K\"ahler geometry, see \cite{BC1}. 
\bt Let $(M,g,\n )$ be a simply connected Hessian manifold 
such that  $\n$ preserves the metric volume form. Then there exists 
a realisation of $(M,g,\n )$ as an improper affine hypersphere 
$\varphi : M \ra \bR^{n+1}$, unique up to 
unimodular affine transformations. Moreover, any automorphism of 
$(M,g,\n )$ has a unique extension to a unimodular affine
transformation of $\bR^{n+1}$ preserving $\varphi (M)$. 
\et 

\pf By the fundamental theorem of affine geometry  
a simply connected pseudo-Riemannian manifold $(M,g,\n )$ with a 
torsion-free connection $\n$ admits a 
Blaschke immersion $\varphi : M \ra \bR^{n+1}$
as a hypersurface with Blaschke metric $g$ and induced connection
$\n$ if and only if the conjugate connection $\bar{\n}$ is 
torsion-free and projectively flat, and if the metric volume form
is $\n$-parallel, see \cite{DNV}.    Moreover,
the Blaschke immersion is unique up to unimodular affine transformations and 
is an improper affine hypersphere if and only if the connection 
$\nabla$ is flat. The assumptions of the fundamental theorem
are satisfied in virtue of Proposition \ref{conjProp} and $\n$ is flat
by Definition \ref{mainDef}. If $\psi : M \ra M$
is an automorphism, then, due to the unicity statement in the fundamental
theorem, there exists a unimodular affine transformation $\alpha : \bR^{n+1}
\ra \bR^{n+1}$ such that $\alpha \circ \varphi = \varphi \circ \psi$.  
It is unique since any affine transformation which fixes an affine frame
is the identity and an affine frame in $\bR^{n+1}$ is determined by a frame in
$T_{\varphi (p)}M$ and the affine normal which is invariant under unimodular
affine tranformations. 
\qed 

\bc  
If $G=Aut(M,g,\n )$ acts transitively on a simply connected 
Hessian manifold $(M,g,\n )$ and $\n$ preserves the metric
volume form, then the Blaschke immersion $\varphi : (M,g,\n )\ra \bR^{n+1}$ 
is a covering
$M=G/H\ra \varphi (M) =Gx_0 =G/G_{x_0}$ over the orbit $Gx_0$ of
a point $x_0\in \varphi (M)$ with respect to an affine action of
$G$ on $\bR^{n+1}$. Moreover, $\varphi (M)\subset \bR^{n+1}$
is an improper affine hypersphere. 
\ec 

\bc Let $(M,g,\n )$ be a Hessian manifold with complete and
positive definite metric $g$ and such that $\n$ preserves the metric volume
form. Then  $g$ is flat and $D=\n$. In particular, any homogeneous
Hessian manifold with positive definite metric and volume 
preserving connection $\n$ is finitely covered
by the product of a flat torus and a Euclidian space. 
\ec  

\pf By the previous theorem, the universal covering of $(M,g,\n )$
can be realised as an improper affine hypersphere with complete and
positive definite Blaschke metric. By the Calabi-Pogerelov theorem, 
see \cite{NS} and references therein, such a hypersurface is a 
paraboloid and the Blaschke metric is flat. Now the existence of 
the finite covering follows from
Bieberbach's theorem. 
\qed 

\section{Geometric structures on the tangent bundle} 
Now we show that the tangent bundle $\pi : N=TM\ra M$ 
of a Hessian
(pseudo-) Riemannian manifold $(M,g,\n )$ has a natural 
(pseudo-) K\"ahler 
structure and the tangent bundle of an affine special real manifold 
has a natural 
special (pseudo-) K\"ahler 
structure. We recall that an {\cmssl (affine) special
(pseudo-) K\"ahler structure} 
$(g,J,\n )$ on a manifold $N$ is given by a (pseudo-) K\"ahler structure
$(g,J)$ and a flat torsion-free symplectic ($\n \omega = \n g \circ J = 0$) 
connection such that
$\n J$ is a symmetric (1,2)-tensor, see \cite{C} and references therein. 

Let $TN= T^hN\oplus T^vN$ be the decomposition of the tangent bundle 
of $N=TM$ into vertical and horizontal
subbundles with respect to the flat connection $\n$ on the Hessian
manifold $(M,g,\n )$. We have a canonical isomorphism
$$T_\xi N = T^h_\xi \oplus T^v_\xi \cong T_{\pi (\xi )}M\oplus 
T_{\pi (\xi )}M.$$ 
Local affine coordinates $x^i$ on $M$ induce
canonical coordinates $(x^i,u^i)$ on $N$ such that any vector $\xi \in TM$ 
is written as $\xi = 
u^{i}\frac{\partial}{\partial x^i}$. Then 
$\partial_i:=\frac{\partial}{\partial x^i}$, 
$\partial_{i'}:=\frac{\partial}{\partial u^i}$
forms a local frame of $T^h$ and $T^v$, respectively. For a 
vector field $X=X^i\frac{\partial}{\partial x^i}$ on $M$ 
we denote by $X^h=X^i\frac{\partial}{\partial x^i}$, 
$X^v=X^i\frac{\partial}{\partial u^i}$ the 
horizontal and vertical lifts of $X$, respectively. 
Then we have the formulas:
$$ [X^h,Y^h]=[X,Y]^h,\quad [X^v,Y^v]=0,\quad [X^h,Y^v]= (\n_XY)^v.$$

The canonical isomorphism 
$\id = \frac{\partial}{\partial u^{i}}\otimes dx^i : T^h_\xi \cong   T_{\pi (\xi )}M \cong T^v_\xi$ defines the
complex structure 
$$J:= \left( \begin{array}{cc} 0& -\id \\
\id & 0\end{array}\right) ,$$
which is integrable since it has constant coefficients
in the coordinates $(x^i,u^i)$. 
We define a natural extension of the metric by
$$g^N:= \left( \begin{array}{cc} g& 0 \\
0& g\end{array}\right) .$$
\bp (cf.\ \cite{S0}) \label{Prop8} For any Hessian (pseudo-) Riemannian 
manifold $(M,g,\n )$ the pair $(g^N,J)$ is a 
(pseudo-) K\"ahler structure on $N=TM$. 
\ep 

\pf It is sufficient to check that the 2-form
$$ \omega = g^N \circ J =  \left( \begin{array}{cc} 0& -g\\
g& 0\end{array}\right) =- g_{ik}(x)dx^i\wedge du^{k} 
$$
is closed. Indeed,
$$d \omega = \frac{\partial g_{ik}}{\partial x^j}dx^i\wedge dx^j\wedge du^{k}
= 0$$
due to the total symmetry of $\n g = \partial g$. \qed 

\bp With respect to the coordinates 
$x^I=(x^i,u^{i'})$, the Christoffel symbols 
$\G_I= (\G_{IJ}^K)$ of the metric $g_N$ are given by 
$$\G_i = \left( \begin{array}{cc} \hat{S}_i& 0\\
0& \hat{S}_i\end{array}\right),\quad \G_{i'} = \left( \begin{array}{cc} 0 & 
-\hat{S}_{i}\\
\hat{S}_{i}& 0\end{array}\right) = J\G_i.$$ 
If $(M,g,\n )$ is special real, that is $\n S=0$, then 
the curvature $R^N$ of the Levi Civita connection $D^N$ is given by
\begin{eqnarray}\label{curvEqu} R^N(X^h,Y^h)&=& R^N(X^v,Y^v)\nonumber \\
&=&\left( \begin{array}{cc} R^D(X,Y)& 0\\
0&R^D(X,Y)\end{array}\right) = -\left( \begin{array}{cc} [\hat{S}_X,\hat{S}_Y]& 0\\
0&[\hat{S}_X,\hat{S}_Y]\end{array}\right)\\ 
R^N(X^h,Y^v) &=& \left( \begin{array}{cc} 0 & {\{ \hat{S}_X,\hat{S}_Y\}}\\ 
-{\{ \hat{S}_X,\hat{S}_Y\}} & 0\end{array}\right) .
\end{eqnarray} 
\ep 

\bc \label{RicCor} Let $(M,g,\n )$ be an affine special real manifold. Then the 
Ricci curvature of the (pseudo-) K\"ahler manifold $(N,g_N,J)$ is
given by
$$ric (X^h+Y^v, X^h+Y^v)= 2 tr \hat{S}_X^2 + 2 tr \hat{S}_Y^2.$$
If the metric $g$ is Riemannian then the Riemannian metric $g_N$
has nonnegative Ricci curvature. The Ricci  curvature is 
strictly positive if and only if  the map 
$X\mapsto S(X,\cdot ,\cdot )$ has trivial kernel.  
\ec 

\pf Since the Ricci form of any K\"ahler manifold is 
given by $\rho (X,Y) = ric (JX,Y)= \frac{1}{2}tr JR(X,Y)$ the formulas
for the Ricci curvature follow from the previous proposition. 
Since $\hat{S}_X$ is symmetric with respect to $g$, 
$tr \hat{S}_X^2\ge 0$ if $g$ is definite.
\qed 

Now we extend the (1,2)-tensor field $\hat{S}$ on $M$ 
considered as tensor on $T^h_\xi N\cong T_{\pi (\xi ) }M$ 
to a (1,2)-tensor $\hat{S}^N$ on $N$ 
such that 
$$\hat{S}^N_J=\hat{S}^NJ=-J\hat{S}^N.$$ 
Then it is given by
\begin{equation} \label{ShatEqu}\hat{S}^N_{X^h} = \left( \begin{array}{cc} \hat{S}_X& 0\\
0&-\hat{S}_X\end{array}\right) ,\quad \hat{S}^N_{X^v} = 
\left( \begin{array}{cc} 0& -\hat{S}_X\\
-\hat{S}_X&0\end{array}\right) .
\end{equation} 
We define a connection $\n^N$ on $N$ by 
$$ \n^N = D^N - \hat{S}^N.$$

\bl \label{connLemma}Let $(N,g^N,J)$ be the K\"ahler manifold associated to
a Hessian manifold $(M,g,\n )$. Then the above connection
$\n^N$ has the following Christoffel symbols with respect to
the coordinates $x^I=(x^i,u^{i'})$:
$$\G_i =  \left( \begin{array}{cc} 0& 0\\
0& 2\hat{S}_i\end{array}\right),\quad \G_{i'} = \left( \begin{array}{cc} 0 & 
0\\
2\hat{S}_i&0\end{array}\right).$$ 

\el 

\bt \label{rmapThm} Let $(M,g,\n = D - \hat{S})$ be an affine special real manifold.
Then $(N,g^N,J,\n^N = D^N - \hat{S}^N)$ is an affine special 
(pseudo-) K\"ahler manifold. 
\et 

\pf Due to Proposition \ref{Prop8} it suffices to show that
the connection $\n^N$ is a) torsion-free, b) symplectic, 
c) special (i.e.\ $\n^NJ$ is symmetric) and d) flat. The properties
a), b) and c) are valid for any Hessian manifold $(M,g,\n )$. 
Indeed the Levi Civita connection $D^N$ is 
torsion-free and preserves $\omega$ and $J$. Therefore:\\  
a) follows from the symmetry $\hat{S}_X^NY=\hat{S}_Y^NX$ for
$X,Y \in TN$.\\
b) follows from the symmetry $\omega (\hat{S}_X^NY,Z) = 
\omega (\hat{S}_X^NZ,Y)$ for all $X,Y, Z \in TN$.\\ 
c) follows from the symmetry of $\hat{S}_X^NY$, since 
$\n^NJ = [\hat{S}_X^N,J]=-2J\hat{S}_X^N$.\\
Now we check that the curvature 
$$R_{IJ}=\partial_I\G_J -
\partial_J\G_I + [\G_I,\G_J]$$
of $\n^N$ vanishes. Using the formula \re{nablaSEqu} for $\n S=0$ 
and the previous lemma we get:
$$\partial_i\G_j =  \left( \begin{array}{cc} 0& 0\\
0& -4\hat{S}_i\hat{S}_j\end{array}\right),\quad \partial_i\G_{j'} = 
\left( \begin{array}{cc} 0 & 
0\\
-4\hat{S}_i\hat{S}_j&0\end{array}\right),\quad \partial_{i'}\G_J =
\frac{\partial}{\partial u^i}\G_J=0.$$ 
Now it is easy to check that $R_{IJ}=0$.\qed 

The same calculation shows:
\bc \label{curvCor} The curvature of the connection $\n^N$ of the 
K\"ahler manifold $(N,g^N,J)$ associated to a Hessian
manifold $(M,g,\n )$ is given by:
$$R(X^h,Y^h)= \left( \begin{array}{cc} 0& 0\\
0& P_{X,Y}-P_{Y,X}\end{array}\right),\quad R(X^v,Y^v)=0,\quad R(X^h,Y^v)=
\left( \begin{array}{cc} 0& 0\\
P_{X,Y}& 0\end{array}\right) ,$$
where $P_{X,Y}Z= g^{-1}(\n_XS)(Y,Z,\cdot ).$
The Ricci curvature  of $\n^N$ is given by:
$$ric (X^h,Y^h)=-tr P_{X,Y},\quad ric(X^v,Y^v)= ric (X^h,Y^v)= 0.$$
\ec 

\bc \label{8.5Cor} Let $(g^N,J,\n^N)$ be the geometric structures on $N=TM$ 
associated to a Hessian manifold $(M,g,\n )$. Then the following
are equivalent:
\begin{enumerate}
\item[(i)] $\n^N$ is flat.
\item[(ii)] $(N,g^N,\n^N)$ is Hessian.
\item[(iii)] $(M,g,\n )$ is special real.   
\end{enumerate}
\ec 

\pf The equivalence of (i) and (iii) follows from the 
previous corollary and it is clear that (ii) implies (i). 
It remains to check that (i) implies (ii). The complete
symmetry of  the tensor field $\n^Ng^N$ follows from the 
symmetry $\hat{S}_X^NY=\hat{S}_Y^NX$, since
\begin{eqnarray*} 
(\n^N_Xg^N)(Y,Z) &=& -(\hat{S}_X^N\cdot g^N)(Y,Z)=g^N(\hat{S}_X^NY,Z)
+g^N(Y,\hat{S}_X^NZ)\\
&=&g^N(\hat{S}_X^NY,Z)-\o^N(J^NY,\hat{S}_X^NZ)\\
&=&g^N(\hat{S}_X^NY,Z)-\o^N(J^N\hat{S}_X^NY,\hat{S}_X^NZ)
= 2g^N(\hat{S}_X^NY,Z).
\end{eqnarray*}
\qed 
\bc \label{Cor9} Under the assumptions of Theorem \ref{rmapThm} the exterior covariant
derivatives of the endomorphism valued one-form $\hat{S}^N$ are given by:
\begin{eqnarray}
d^{\n^N}\hat{S}^N(X^h,Y^h) &=& d^{\n^N}\hat{S}^N(X^v,Y^v)=
-2\left( \begin{array}{cc} [\hat{S}_X,\hat{S}_Y] & 
0\nonumber \\
0&[\hat{S}_X,\hat{S}_Y]\end{array}\right)\
\nonumber \\
d^{\n^N}\hat{S}^N(X^h,Y^v) &=& 2\left( \begin{array}{cc} 0 & 
\{ \hat{S}_X,\hat{S}_Y\} \nonumber \\
-\{ \hat{S}_X,\hat{S}_Y\}&0\end{array}\right)\nonumber \\
d^{D^N}\hat{S}^N=0
\end{eqnarray}
\ec 

\pf The proof follows from \re{curvEqu}, $R^{\n^N}=0$ and the formulas 
\begin{eqnarray}
d^{\n^N}\hat{S}^N &=& R^{N}-R^{\n^N}-[\hat{S}^N,\hat{S}^N]\nonumber \\
d^{D^N}\hat{S}^N &=& R^{N}-R^{\n^N} + [\hat{S}^N,\hat{S}^N] \nonumber. 
\end{eqnarray}
\qed 

\bc \label{Cor10}
Under the assumptions of Theorem \ref{rmapThm} the curvature of the
Levi-Civita connection of $g_N$ is given by 
\[ R^N=-[\hat{S}^N,\hat{S}^N].\]
\ec 

\bd The map which to any affine special real manifold $(M,g,\n )$
associates the affine special K\"ahler manifold $\mathbf{r}(M):= (N,g^N,J,\n^N )$ 
is called
the {\cmssl (affine) r-map}.  
\ed 
\bc  Let $(M,g,\n )$ be a special real manifold which locally 
admits a homogeneous Hesse potential $h=\frac{1}{6}S_{ijk}x^ix^jx^k$.  Then the 
K\"ahler manifold $\mathbf{r}(M)$ is not flat. 
\ec   

\pf By \re{curvEqu} $R^N=0$ is equivalent to $\hat{S}_X\hat{S}_Y=0$ 
for all $X,Y\in TM$. This is impossible since, by Proposition 
\ref{PropAffSpecMan}, 
locally we can identify the Hessian manifold with a domain $(V(x_0),
g,\partial )$ and 
$2g\circ \hat{S}_X=S(X,\cdot ,\cdot ) = \partial^2 h (X)= g_X$ 
is invertible for all $X \in V(x_0)$. 
\qed 

Flat special K\"ahler manifolds were classified in \cite{BC2}. By the
corollary, they cannot be obtained from a homogeneous cubic
polynomial by the r-map. 

\bt \label{characThm}
Let $(M,g,\n )$ be a Hessian manifold of dimension $n$ and\linebreak  
$(N=TM,g^N,J,\n^N)$ the corresponding (pseudo-) K\"ahler manifold
with the connection defined above. Then  
\begin{itemize}
\item[(i)] The decomposition $TN=T^v+T^h=T^v+JT^v$ is a decomposition
into two orthogonal 
integrable Lagrangian distributions which are totally geodesic
and flat with respect to $\n^N$. The horizontal distribution is also totally
geodesic with respect to $D^N$. 
\item[(ii)] For any leaf $L=M(\xi )$, $\xi \in N$, of the horizontal
distribution there exists an involution $\s=\s_L\in {\rm Aut}(N,g^N,J,\n^N)$
which preserves the vertical and horizontal foliations and 
such that $L=N^{\s }$. 
\item[(iii)] The group generated by products $\s_L\circ \s_{L'}$ preserves
each fiber $T_pM$, $p\in M$, and acts as the translation group of the 
fiber $T_pM$.   
\end{itemize}
\et 

\pf (i) follows from the formulas for the Christoffel symbols and curvature 
of $D^N$ and $\n^N$ and from the formulas for 
$g^N=g_{ij}(x)(dx^idx^j+du^idu^j)$,
$J$ and $\omega =g^N\circ J$. In each coordinate domain one 
can check directly that the reflection 
$\s_U: (x,u) \mapsto (x,-u+2u_0)$ with respect to an open domain 
$U=\{ u=u_0\}\subset L$ in a leaf $L$ 
has the  properties claimed in (ii) . The reflections $\s_U$ coincide on 
overlaps and, hence, define the global reflection $\s_L$. (iii) follows from
the fact that the product of two central symmetries in the affine space
$T_pM$ is a parallel translation. \qed 

The converse can be stated as follows. 
\bt \label{characconvThm}
Let $(N,J,g^N)$ be a 2n-dimensional pseudo-K\"ahler manifold which admits 
a free holomorphic and isometric action of the vector group
$\bR^n$ with Lagrangian orbits such that the projection $\pi :N\ra 
N/\bR^n=M$ is a  trivial (principal) bundle. Then there exists an induced 
pseudo-Riemannian metric $g$ and flat connection $\n$ on $M$ such that
$(M,g,\n )$ is Hessian and $N$ is identified with $TM$ with the 
pseudo-K\"ahler structure  $(J,g^N)$ induced from $(M,g,\n )$ by Proposition
\ref{Prop8}. If moreover $(N,J,g^N,\n^N)$ is special K\"ahler and the 
Killing vector fields $U_i$ of the above action are $\n^N$-parallel
along the Lagrangian orbits, then 
$(M=N/\bR^n,g,\n )$ is special real and $(N,J,g^N,\n^N)$ is obtained
from $(M=N/\bR^n,g,\n )$ by the r-map. 
\et 

\pf We denote by $U_i$ the commuting vector fields on $N$ which are the 
generators
of the action of $\bR^n$. The holomorphicity of the $U_i$ implies that
the vector fields $U_i,X_j=-JU_j$ commute. They are linearly independent
since the distribution $T^vN:= {\rm span}\{ U_i| i=1,\cdots ,n\}$ is 
Lagrangian. There exist local coordinates $(x^i,u^j)$ such that
$U_i = \frac{\partial }{\partial u^i}$ and $X_j = 
\frac{\partial }{\partial x^j}$. In these coordinates
$$g^N=g_{ij}(dx^idx^j+du^idu^j),$$
where the functions $g_{ij}=g_{ij}(x)$ depend only on the $x^i$, since
the $U_i$ are Killing vector fields. Since $U_i$ is a 
holomorphic Killing vector field, it is also 
symplectic with respect to the K\"ahler form. 
This implies that $X_i=-JU_i$ is a gradient vector field, hence
the one-form $g(X_i,\cdot )$ is closed and 
$\frac{\partial }{\partial x^i}g_{jk}$ is completely symmetric. 
This shows that $(M,g,\n )$ is a Hessian manifold, where $g=(g_{ij}(x))$ 
is the metric on $M$ which makes $N\ra M$ a pseudo-Riemannian submersion and 
$\n$ is the flat connection induced by the flat connection
on the leaves of the distribution $T^h:=JT^v$ with parallel vector fields  
$\frac{\partial }{\partial x^i}$. Identifying $M$ with a 
section of the trivial bundle $N\ra M$ we can identify $N$ with $TM$.  
It is clear that the 
pseudo-K\"ahler structure on $N=TM$ is obtained by Proposition
\ref{Prop8} from the Hessian manifold $(M,g,\n )$.

In the special K\"ahler case the assumption $\n^N_{U_i}U_j=0$ implies
$$\n_{JU_i}^N(JU_j)=(\n_{JU_i}^NJ)U_j + J\n_{JU_i}^NU_j=
-J(\n_{U_j}^NJ)U_i+J(\n_{U_j}^NJ)U_i=0.$$
Using the fact that for a special K\"ahler manifold $\hat{S}^N:=
D^N-\n^N=-\frac{1}{2}J\n^NJ$, one can easily check  
that $\n^N$ is the connection from Lemma \ref{connLemma}. 
Now Corollary \ref{curvCor} shows that $\n S=0$, i.e.\ $(M,g,\n )$ 
is special real. \qed  
As an application we prove the following non-existence result.
\bt \label{lastThm} There is no compact simply connected special
real manifold of positive dimension.
\et 

\pf Let $(M,g,\n)$ be a compact simply connected special
real manifold. We first prove that the metric of the special
K\"ahler manifold $(N=TM,g^N,J,\n^N)$ is complete. 
As in the proof of Corollary \ref{immCor}, 
we have a $\n$ parallel coframe $\xi^i=dx^i$, $i=1,\ldots ,n$. We denote
by $\partial_i$ the dual global frame on $M$ and by $u^i$
the corresponding globally defined functions on $N=TM$, which are the 
coordinates of a vector with respect to the frame $\partial_i$. 
There exists a free isometric action of $\bR^n$ on $N=TM$
given by $u^i \mapsto u^i+c^i$, whose orbits are the tangent spaces. 
The quotient of $(N,g^N)$ by the lattice of integral translations
is a compact, hence, complete Riemannian manifold. Thus  
$(N,g^N)$ is complete as the universal covering of a 
complete Riemannian manifold. 

By \cite{BC1}, $(N,g^N,\n^N)$ can be realised as an improper affine hypersphere
with complete Riemannian metric. In fact, one can check that
the integrability conditions\footnote{These are: $\n^N$ is flat, torsion-free 
and preserves the metric volume and the conjugate connection
$\bar{\n}$ is (projectively) flat and torsion-free.} for the existence of
an affine hypersphere immersion $\varphi : N \ra \bR^{2n+1}$
follow from the properties of a special K\"ahler manifold.  By the 
Calabi-Pogorelov theorem, see 
\cite{NS} and references therein, the metric $g^N$ admits a quadratic Hesse
potential and the Levi-Civita connection of $g^N$ coincides with $\n^N$. 
Using the formulas \re{ShatEqu}, 
this implies that the Levi-Civita connection of the special real
manifold coincides with $\n$ and, hence, is flat. Since there is no
compact simply connected Riemannian manifold of positive dimension
we obtain the theorem. \qed

\end{document}